
\def\date{\ifcase\month\or January\or February \or March\or April\or May %
\or June\or July\or August\or September\or October\or November           %
\or December\fi\space\number\day, \number\year}                          %
\def\date{08-03-2012}      
\magnification=\magstep1   
\hsize=16truecm            
\vsize=24truecm            
\parindent=0pt             
 \line{\hfill\date}        
\catcode`\@=11             
\newfam\msbfam             
\newfam\euffam             
\font\tenmsb=msbm10        
\font\sevenmsb=msbm7       
\font\teneuf=eufm10        
\font\seveneuf=eufm7       
\textfont\msbfam=\tenmsb   
\scriptfont\msbfam=\sevenmsb%
\textfont\euffam=\teneuf   
\scriptfont\euffam=
\seveneuf                  
\def\I{{\fam\msbfam I}}    
\def\L{{\fam\euffam L}}    
\def\N{{\fam\msbfam N}}    
\def\T{{\fam\msbfam T}}    
\def\R{{\fam\msbfam R}}    
\def\Z{{\fam\msbfam Z}}    
\def\g{{\fam\euffam g}}    
\catcode`\@=12             
\def\.{{\cdot}}            
\def\<{\langle}            
\def\>{\rangle}            
\def\({\big(}              
\def\){\big)}              
\def\hat{\widehat}         
\def\fin #1{#1_{\rm fin}}  
\def\implies{
\hbox{$\Rightarrow$}}      
\def\im{\mathop{\rm im}    
\nolimits}                 
\def\id{\mathop{\rm id}    
\nolimits}                 
\def\Aut{\mathop{\rm Aut}  
\nolimits}                 
\def\Hom{\mathop{\rm Hom}  
\nolimits}                 
\def\pr{\mathop{\rm pr}    
\nolimits}                 
\def\defi{\buildrel\rm def 
\over=}                    
\def\ssk{\smallskip}       
\def\msk{\medskip}         
\def\bsk{\bigskip}         
\def\nin{\noindent}        
\def\cen{\centerline}      
\font\smc=cmcsc10          
\def\lead{\leaders\hbox to 1.5ex{\hss${.}$\hss}\hfill}
\def\arr{\hbox to 40pt{\rightarrowfill}}
\def\larr{\hbox to 40pt{\leftarrowfill}}

\newcount\litter                                         
\def\newlitter#1.{\advance\litter by 1                   
\edef#1{\number\litter}}                                 
\def\sdir#1{\hbox{$\mathrel\times{\hskip -4.6pt          
 {\vrule height 4.7 pt depth .5 pt}}\hskip 2pt_{#1}$}}   
\def\qedbox{\hbox{$\rlap{$\sqcap$}\sqcup$}}              
\def\qed{\nobreak\hfill\penalty250 \hbox{}               
\nobreak\hfill\qedbox\vskip1\baselineskip\rm}            

\long\def\alert#1{\parindent2em\smallskip\line{\hskip\parindent\vrule%
\vbox{\advance\hsize-2\parindent\hrule\smallskip\parindent.4\parindent%
\narrower\noindent#1\smallskip\hrule}\vrule\hfill}\smallskip\parindent0pt}

\newlitter\bour.
\newlitter\bred.
\newlitter\diac.
\newlitter\ert.
\newlitter\EK.
\newlitter\ERE.
\newlitter\ER.
\newlitter\gust.
\newlitter\hew.
\newlitter\HM.
\newlitter\HMpro.
\newlitter\hofmosii.
\newlitter\lesc.
\newlitter\levai.
\newlitter\neum.
\newlitter\neumtwo.
\newlitter\shalev.
\newlitter\tD.

\cen{\bf  The probability that x and y commute} 
\msk 
\cen{\bf in a compact group} 
\msk 
\cen{Karl H. Hofmann and Francesco G. Russo}

\bsk {\smc Abstract.}
We show that a compact group $G$ has finite conjugacy classes, i.e.,  
is an FC-group if and only if its center
$Z(G)$ is open if and only if its commutator subgroup $G'$ is finite.
Let $d(G)$ denote the Haar measure of the set of all
pairs $(x,y)$ in $G\times G$ for which $[x,y]=1$; this, formally,
 is the probability that two randomly picked elements commute. 
We prove that $d(G)$ is always
rational and that it is positive if and only
if  $G$ is an extension of an FC-group by a finite group. This
entails that $G$ is abelian by finite.  The proofs 
involve measure theory, transformation groups, 
Lie theory of arbitrary compact
groups, and representation theory of compact groups.
  Examples and  references to the
history of the discussion are given at the end of the paper.

\msk {\smc MSC 2010}: Primary 20C05, 20P05; Secondary 43A05.

\msk {\smc Key words}: Probability of commuting pairs, commutativity
degree, FC-groups,  compact groups, Haar measure.

\bsk

\cen{\smc 1. Introducing the problem}
\msk
We fix a compact group $G$ and 
 let $\nu$ denote Haar measure of $G$.

Question. {\it What is the probability that a randomly picked pair
$(x,y)\in G\times G$ of elements of $G$ has the property that $x$
and $y$ commute?} 

\ssk Technically speaking: Derive useful information from possible
knowledge on the real number
$$(\nu\times\nu)(\{(x,y)\in G\times G: [x,y]=1\}).$$

We shall develop  a
theory which will show that in a compact group $G$ 
the probability of the commuting of two randomly picked 
elements $x$ and $y$
is positive only in very special groups,
namely, those that have a certain characteristic 
open abelian normal subgroup. Our theory will allow us to 
compute this probability explicitly and to establish that it
is always a rational number.

\msk

The FC-{\it center} of a group is the set of elements whose conjugacy class is
finite; a group is an FC-{\it group} if it agrees with its FC-center. 
We shall prove the following Theorems:

\msk
\bf
Theorem 1.1.~\rm(The structure of compact FC-groups)\quad\it 
A compact group $G$ is an {\rm FC}-group if and only
if its center is open, that is, $G$ is central by finite, if and only
if its commutator subgroup is finite, that is, $G$ is finite by abelian.  
\msk
\bf
Theorem 1.2.~\rm(The structure of compact groups with frequent commuting of 
elements)\quad \it 
Let $G$ denote a compact group with Haar measure
$\nu$ and {\rm FC}-center $F$.
Then the following conditions are equivalent:
{\parindent2em

\item{$(1)$} $(\nu\times\nu)(\{(x,y)\in G\times G:[x,y]=1\})>0$.

\item{$(2)$} $F$ is open in $G$.

\item{$(3)$} The center $Z(F)$ of $F$ is open in $G$.

}\rm\bsk

In our proofs,  we shall
have to draw from different branches of group theory such as
the theory of compact groups, integration on compact spaces and groups, 
transformation groups, representation theory, and FC-groups.

\bsk

\cen{\smc 2. Actions and product measures}

\msk
It is conceptually simpler to start by considering a more general situation.
Let $(g,x)\mapsto g\.x:G\times X\to X$ be a continuous action
$\alpha$ of a compact group $G$ on a compact space $X$.
All spaces in sight are assumed to be Hausdorff.
We specify a Borel probability measure $P$ on $G\times X$ and
 discuss the probability that a group element $g\in G$ fixes
a phase space element $x\in X$ for a pair $(g,x)$, randomly picked
from $G\times X$, that is, that  $g\.x=x$.
We define
$$E\defi \{(g,x)\in G\times X: g\.x=x\},$$
that is, $E$ is the equalizer of the two functions
$\alpha,\pr_X:G\times X\to X$ and is therefore a closed
subset of $G\times X$.
\bsk
Let $G_x=\{g\in G:g\.x=x\}$ be the isotropy (or stability) group at $x$
and let $X_g=\{x\in X: g\.x=x\}$ be the set of points fixed under the action
of $g$.
We note that $G_{g\.x}=gG_xg^{-1}$. The function $g\mapsto g\.x:G\to G\.x$
induces a continuous equivariant bijection $G/G_x\to G\.x$ which,
due to the compactness of $G$, is a homeomorphism.
 We have

$$\eqalign{%
E&=\{(g,x):g\in G_x,\, x\in X\}=\bigcup_{x\in X}G_x\times\{x\}\cr
&=\{(g,x): g\in G,\, x\in X_g\}=\bigcup_{g\in G}\{g\}\times X_g
\subseteq G\times X\.\cr}$$

\bsk

Now we assume that  $\mu$ and $\nu$ are Borel probability  measures on
$G$ and $X$, respectively, and that $P=\mu\times\nu$
is the product measure. For information on measure theory
the reader may refer to [\bour]. 
Let $\chi_E\colon G\times X\to \R$ be the characteristic function
of $E$.
We define the function $m\colon X\to \R$, $m(x)=\mu(G_x)$.
Then by the Theorem of Fubini we compute

$$\leqalignno{P(E)&=\int_{G\times X}\chi_E(g,x)dP
=\int_X\left(\int_G\chi_E(g,x)d\mu(g)\right)d\nu(x)&(*)\cr
&=\int_X\mu(G_x)d\nu(x)=\int m\,d\nu.\cr} $$
Likewise
$$\leqalignno{%
P(E) &= \int_{G\times X}\chi_E(g,x) dP
      =\int_G\left(\int_X\chi_E(g,x)d\nu(x)\right)d\mu(g)&(**)\cr
     &=\int_G \nu(X_g) d\mu(g).\cr}$$

We now see from $(*)$ that $P(E)>0$ implies, 
firstly, that there is at least one
$x\in X$ such that
$\mu(G_x)=m(x)>0$ holds, and that, secondly, the set of all $x$ for which 
$m(x)>0$ has positive $\nu$-measure. Likewise, there is at least one 
$g\in G$ such that $\nu(X_g)>0$ and that the set of all of these $g$
has positive $\mu$-measure.
At this point we introduce a terminology which we shall retain and
use 

\msk
\bf
Definition 2.1.\quad\rm  Let $\cal C$ be a set of subgroups of a compact
group $G$, such as the set of all closed subgroups 
or all subgroups whose underlying
set is a Borel subset of $G$. We shall say that
a Borel probability measure  $\sigma$ 
on $G$ {\it respects $\cal C$-subgroups}
if every subgroup $H\in{\cal C}$  with
$\sigma(H)>0$ is open.
\msk
Recall that an open subgroup $H$ of a topological group $G$,
being the complement of
all the cosets $gH$ for $g\notin H$, is closed and that
it contains the {\it identity component} $G_0$ of $G$. If $G$ is
compact, then $H$ has finite index in $G$.

\msk

We claim that 

\centerline{ Haar measure $\mu$ on a compact group $G$ respects Borel subgroups.}

Indeed, assume
that $H$ is a Borel subgroup of $G$ with $\mu(H)>0$.
Then by [\hew], p.~296, Corollary 20.17, $H=HH$ contains a nonvoid open
set and thus is open.
\ssk 

\msk

If, in the cirumstances discussed here, $\mu$ respects closed subgroups,
and $\mu(G_x)>0$ then  the subgroup
$G_x$ is open in $G$ hence contains $G_0$. 
If $G$ is a compact Lie group, then  the condition
$G_0\subseteq G_x$ is  also sufficient for the openness of the subgroup 
$G_x$  of $G$,
 since the identity component 
of a Lie group is open.   
If $G_x$ is open,  then $G\.x\cong G/G_x$ is discrete
and compact, hence finite. 

\ssk
Let now
$$F\defi\{x\in X: |G\. x|<\infty \}
=\{x\in X: G_x\hbox{ is open}\}.\leqno(\dag)$$
These remarks require that henceforth whenever the measure $\mu$ occurs
 we shall assume that 
$\mu$ respects closed subgroups.
\msk

\bf Lemma 2.2.\quad \it Let $G$ be a compact group acting on a compact
space $X$.
\ssk
{\parindent=2em%
\item{\rm(i)}  If $G$ is a Lie group, then for each $x\in X$ there 
          is an open invariant
          neighborhood $U_x$ of $x$ such that all isotropy groups
          of elements in $U$ are conjugate to a subgroup of 
          the isotropy subgroup $G_x$.
\ssk\item%
{\rm(ii)} Under these circumstances, $m$ takes its
          maximum on $U_x$ in  $x$. That is, \hfill\break 
          $m^{-1}(\left]-\infty,m(x)\right])$ is a
          neighborhood of $x$. 
          In particular, $m$ is upper semicontinuous.

\ssk\item%
{\rm(iii)} If $G$ is a Lie group, then the subspace $F$ of $X$ is compact.

\ssk\item%
{\rm(iv)} If $G$ is an arbitrary compact group, then
the subspace $F$ of $X$ is an {\rm F}$_\sigma$, that is, a countable 
union of closed subsets and thus is a Borel subset.

\ssk\item%
{\rm(v)} There is a compact commutative monoid $X$ with compact group $G$ of 
invertible elements such that for the action of $G$ on $X$ by multiplication,
the set $F=X\setminus G$ is a nonclosed $F_\sigma$. 
}
\ssk
\bf Proof. \quad\rm  Assertion (i) is a consequence of the Tube Existence
Theorem 
(see e.g. [\bred], p.~86, Theorem 5.4, or [\tD], p.~40, Theorem 5.7).

\ssk
(ii)   Immediate from (i) and from $m(u)=\mu(G_u)$.

\ssk
(iii)  We recall that $y\in F$ if and only if $m(y)>0$. Hence $F$ is
the complement of $m^{-1}(0)$. By conclusion (ii), however,
$m^{-1}(0)$ is open. Thus $F$ is closed and therefore compact.

\ssk
(iv) For a natural number $n\in \N$ we set $F(n)=\{x\in X: |G\.x|\le n\}$.
We claim that $F(n)$ is closed in $X$ for all $n\in \N$. 
Since $F=\bigcup_{n=1}^\infty F(n)$,
this claim will prove assertion (iv).

We prove the claim by contradiction and suppose that there is an $n\in\N$
such that $\overline{F(n)}$ contains an $x'\notin F(n)$. Then there exist
elements $g_1,\dots,g_{n+1}\in G$ such that $|\{g_1\.x',\dots,g_{n+1}\.x'\}|=
n+1$. Now we find a compact normal subgroup $N$ of $G$ such that 
$G/N$ is a Lie group and that $Ng_j\.x'\cap Ng_k\.x'=\emptyset$ for 
all $j\ne k$ in $\{1,\dots,n+1\}$.  
\ssk
The Lie group $G/N$ acts on $X/N=\{N\.x:x\in X\}$ via $(gN){\circ}(N\.x)=N\.(g\.x)$.
By what we just saw $|(G/N){\circ}(N\.x_0)|\ge n+1$. On the other hand,
$F_N(n)=\{N\.x\in X/N: |(G/N){\circ}(N\.x)|\le n\}$ is closed by (ii)
above. Since the orbit map  $\pi_N\colon X\to X/N$ is continuous and
$\pi(F(n))\subseteq F_N(n)$ we have $N\.x_0=\pi(x_0)\subseteq 
\overline{F_N(n)}=F_N(n)$. Thus, by the definition of $F_N(n)$ we
have $|(G/N){\circ}(N\.x_0)|\le n$. This contradiction proves the claim.

\ssk
(v) Let $\I$ denote the unit interval $[0,1]$ under ordinary multiplication,
a compact connected topological monoid. 
Set $\I_0=[0,1-1/p]$, $\I_1=]1-1/p,1-1/p^2]$,
$\dots,\quad\I_n=]1-1/p^n,1-1/p^{n+1}], \dots, \I_\infty=\{1\}$.
We form the compact 
topological product monoid $S=\I\times\Z_p$ for the additive group
$\Z_p$ of $p$-adic integers. 
 For $r\in\I$ set 
$$J_r=\cases{\Z_p, &if    $r\in\I_0$,\cr 
            p^n\Z, &if   $r\in\I_n$, $n=1,2,\dots$,\cr
             \{0\}, &if  $\in\I_\infty$. \cr}$$
The binary relation $R$, whose cosets are
$R(t,z)=\{t\}\times(z+J_t)$ is a closed congruence relation.
Therefore, $X\defi S/R$ is a compact connected abelian monoid 
with zero $R(0,0)$ whose 
group of units is $G\defi(\{1\}\times\Z_p)/R\cong\Z_p$.
For $t<1$ in $\I$ and $x=R(t,z)$  we have 
$G\.x=(\{t\}\times(z+\Z_p/I_p)\cong\Z_p/p^n\Z_p\cong\Z/p^n\Z$ 
for $t\in\I_n$, $n=0,1,\dots\,$.

In particular, considering $X$ as a $G$-space under multiplication, 
the space of
finite orbits $F$ is $([0,1[\times G)/R$ is not closed in $X$. \qed

\ssk
Statement 2.2(v) shows that  2.2(iv)
cannot be improved to read that $F$ is closed.

For $p=2$, the space $X$
 is the standard binary tree with $G$ as the Cantor set of 
leaves. Compact monoids like $X$ above
were considered in [\hofmosii] rather generally 
under the name {\it cylindrical 
semigroups}; for our construction see in particular D-2.3.3ff
on p.~241.
\msk

Recall that we assume that $\mu$ respects closed subgroups. 
Now that we know  that $F\subseteq X$ is a Borel set, hence is 
$\nu$-measurable, we can state
that,  regardless of any particular property of $\nu$,
the function $m$ satisfies 
$$P(E)=\int_Xm d\nu=\int_X\chi_H\.m d\nu=\int_{x\in F}m(x)d\nu(x).$$
Here $x\in F$ implies $0<m(x)=\mu(G_x)=1/|G/G_x|\le 1$.

\msk
\bf

Lemma 2.3.\quad \it
If $\mu$ respects closed subgroups and $P(E)>0$,  
then $0<P(E)\le\nu(F)$.
 In particular, $F\ne\emptyset$.
           
\msk
\bf
Proof.\quad \rm
 We have seen that  $P(E)=\int_Fm\,d\nu$ since $F$ is Borel measurable.
The Lemma then follows from this fact and $m(x)\le1$ for $x\in F$.\qed
\msk
We shall say that the group $G$ {\it acts automorphically} on $X$ if $X$ is a 
compact group and $x\mapsto g\.x:X\to X$ is an automorphism for all 
$g\in G$.

\msk
\bf
Lemma 2.4. \quad\it Assume that $G$ and $X$ are compact groups and
assume the following hypotheses:

{\parindent2em
           \item{\rm(a)} $G$ acts automorphically on $X$.
           \item{\rm(b)} $\mu$ respects closed subgroups.
           \item{\rm(c)} $\nu$ respects Borel subgroups or else
                         $X$ is a Lie group and $\nu$ respects closed
                         subgroups. 
           \item{\rm(d)} $P(E)>0$.
}

Then $F$ is an open, hence closed  subgroup of $X$. 
\msk
\bf
Proof.\quad\rm  Let $x, y\in F$. Then $G\.x$ and $G\.y$ are finite sets
by the definition of $F$. Now $G(xy^{-1})=\{g\.(xy^{-1}):g\in G\}=
\{(g\.x)(g\.y)^{-1}:g\in G\}\subseteq \{(g\.x)(h\.y)^{-1}:g,h\in G\}=
(G\.x)(G\.y)^{-1}$, and the last set is finite as a product 
of two finite sets. Thus  $xy^{-1}\in F$ and $F$ is a subgroup.
By Lemma 2.3, $\nu(F)>0$. Then by Lemma 2.2(iii),(iv) and 
the kind of subgroups respected by  $\nu$, we  
conclude that $F$ is an open subgroup.\qed 
\msk
If $G$ acts automorphically on a compact group $X$, we let 
$\pi\colon G\to \Aut \ X$ be the representation given by 
$\pi(g)(x)= g\.x$. Let $\id_X$ denote the identity function of
$X$. Then the fixed point set
$X_g$ is the equalizer of the morphisms $\pi(g)$ and $\id_G$ and
is therefore a closed subgroup of $X$. Let $I\subseteq G$ denote
the set of all $g\in G$ for which $X_g$ has inner points.

\msk\bf Lemma 2.5.\quad\it
Assume that $G$ and $X$ are compact groups and
assume the following hypotheses:

{\parindent2em
           \item{\rm(a)} $\mu$ and  $\nu$ are the Haar measures on $G$
                         and $X$, respectively.
           \item{\rm(b)} $G$ acts automorphically on $X$. 
           \item{\rm(c)} $G$ is finite.
}

Then $P(E)={1\over|G|}\.\sum_{g\in I} |X/X_g|^{-1}$. 
In particular, $P(E)$ is a rational number.
\msk
\bf
Proof.\quad\rm 
By $(**)$ above and the fact that Haar measure on a finite group $G$
is counting measure with $\mu(\{g\})=|G|^{-1}$, we have
$P(E)={1\over|G|}\.\sum_{g\in G}\nu(X_g)$.
If a closed subgroup $Y$ of the compact group $X$ has no inner points, its
Haar measure $\nu(Y)$ is zero. If it has inner points, it is open
and its measure $\nu(Y)$ is the reciprocal of its index, that is
$\nu(Y)=|X/Y|^{-1}$. Hence $P(E)={1\over|G|}\.\sum_{g\in I} |X/X_g|^{-1}$ 
and the assertion follows.\qed 

\msk

Our conclusions sum up to the following result:
\msk
\bf
Proposition 2.6.\quad\it Let $G$ and $X$ be compact groups and assume that
$G$ acts automorphically on $X$. Let $\mu$ and $\nu$ be normalized 
positive  Borel measures on $G$ and $X$, respectively. Define
$$ E=\{(g,x)\in G\times X: g\.x=x\}\subseteq G\times X,$$
whence $(\mu\times\nu)(E)=\int_{x\in F}\mu(G_x) d\nu(x)$.
Assume that $\mu$ respects closed subgroups and that
 $\nu$ either  respects Borel subgroups of $X$ or else
$X$ is a Lie group and $\nu$ respects closed subgroups,
then the following statements are  equivalent:
{\parindent2em 
\item{$(1)$} $(\mu\times\nu)(E)>0$.

\item{$(2)$} The  subgroup $F\le X$ of all elements
with finite $G$-orbits is open and thus has finite index in $X$.
} \qed \rm
\bsk
\cen{\smc 3. Action by inner automorphism}
\msk

The main application of this general situation will be the case
of a  compact group  $G$ and the automorphic action of $G$ on $X=G$
via inner automorphisms:

$$(g,x)\mapsto g\.x=gxg^{-1}: G\times X\to X.$$

The orbit $G\.x$ of $x$ is the conjugacy class $C(x)$ of $x$,
and the isotropy group $G_x$ of the action at $x$ is the centralizer
$Z(x,G)=\{g\in G: gx=xg\}$ of $x$ in $G$. The set $E$ is the set
$D=\{(x,y)\in G\times G: [x,y]=1\}$,  and $F$ is
the union of all finite conjugacy classes. 
In particular, $F$ is a characteristic  {\rm F}$_\sigma$ 
subgroup of $G$ whose
elements have finite conjugacy classes, that is,
the {\rm FC}-center of $G$. Recall that a group agreeing with 
its FC-center is called an {\rm FC}-group. 

In this setting, Proposition 2.6 has the following consequence:

\msk
\bf 
Corollary 3.1.\quad\it Let $G$ be a compact group and 
let $\mu$ and $\nu$ be  Borel probability
measures on $G$ and assume that $\mu$ respects closed subgroups
and that $\nu$ 
respects Borel subgroups or, if $G$ is a Lie group, that 
$\nu$ respects closed subgroups.

\ssk
  Let $F$ be the {\rm FC}-center
of $G$.  Then $F$ is an {\rm F}$_\sigma$ and we define
$$ D=\{(g,x)\in G\times G: [g,x]=1\}\subseteq G\times G.$$
Then 
$$P(D) =\int_{x\in F}\mu(Z(x,G)) d\nu(x),$$ 
and the following statements are 
equivalent:
{\parindent2em 
\item{$(1)$} $(\mu\times\nu)(D)>0$.

\item{$(2)$} $F$  is open in $G$ 
and thus has finite index in $G$. \qed

}
\rm
Moreover, under these conditions, $Z(F,G)$ contains the identity
component $G_0$, and the profinite group $G/Z(F,G)$ is acting
effectively on $F$ with orbits being exactly the finite
conjugacy classes of $G$. We shall see later in 3.10, 
that the center $Z(F)$ 
is an open subgroup, whence 
$Z(F,G)$, containing $Z(F)$, is an open subgroup of $G$. Therefore
$G/Z(F,G)$ will in fact turn out to be finite.

\bsk
 In a nontechnical spirit, let us call a subgroup $U$ of $G$ 
{\it large} if it is open, hence closed and with finite index
in $G$.
In 3.1 (2) we state that $F$ is a 
large characteristic subgroup which is itself an {\rm FC}-group.
Therefore,  for finding large abelian subgroups of $G$
under the equivalent
hypotheses of 3.1,   it is 
sufficient to concentrate on
compact {\rm FC}-groups which we shall do now; and while we are
classifying compact FC-groups,  measures will not play any role.
\ssk
In a compact {\rm FC}-group $G$, the identity component $G_0$ is abelian
by the structure theory of compact connected groups ([\HM], 9.23ff).
Abbreviate the centralizer $Z(G_0,G)$  of $G_0$ by $C$.
The action by inner automorphisms of $G$ on $G_0$ induces an 
effective automorphic action
of the profinite group $G/C$
on the compact connected abelian group $G_0$ with finite orbits.
Let the Lie algebra
$\g=\L(G_0)$ be defined as $\Hom(\R,G_0)\cong \Hom(\hat{G_0},\R)$ 
(see [\HM], Proposition 7.36, Theorem 7.66).
The Lie algebra $\g$ is a weakly complete vector space, which is isomorphic as
a topological vector space to the complete locally convex space $\R^I$
for a suitable set $I$ whose cardinality is the rank of the torsionfree
abelian group $\hat{G_0}$.  
There is a natural morphism  $\Aut G_0\to \Aut \g$: Indeed each automorphism 
$\alpha$ of $G_0$ induces an automorphism $\L(\alpha)$ of $\g$ as follows:
Let $X\colon \R\to G_0$, $t\mapsto X(t)$, be a member of $\g$, then 
$\L(\alpha)(X)\colon \R\to G_0$ is defined by $\L(\alpha)(X)=\alpha\circ X$.
An element $g\in G$ induces an inner automorphism $I_g$ on $G_0$
via $I_g(x)=gxg^{-1}$ and the function $g\mapsto I_g:G\to\Aut \g$ factors
through the quotient $G\to G/C$ giving us a chain of representations
$$G\to G/C\to \Aut G_0 \to \Aut \g,$$
whose composition yields a continuous representation 
$$\pi\colon G/C\to \Aut \g,\quad \pi(gC)(X)= I_g\circ X, 
\hbox{ that is, } \pi(gC)(X)(t) = gX(t)g^{-1}.$$

We claim that the fact, that every element of $G$ has finitely
many conjugates, implies  for each $X\in \g$ that $\pi(G/C)(X)\subseteq \g$ is 
contained in a finite dimensional vector subspace of  $\g$.
That is, we maintain that 
the $G/C$-module $\g$ (see [\HM], Definition 2.2) satisfies 
$\g_{\rm fin} =\g$. 
(See [\HM], Definition 3.1.) 
\ssk
In order to prove this claim in the next lemma we recall that a subgroup of a 
topological group is called {\it monothetic} if it contains a dense cyclic
subgroup, and {\it solenoidal} if it contains a dense one-parameter subgroup.
\msk
\bf
Lemma 3.2.\quad\it Let $\Gamma$ be a compact group acting automorphically
on a compact  group $G$.
 Assume that all orbits of $\Gamma$
on $G$ are finite. Then the following conclusions hold:
{\parindent2em 

\item{\rm(i)} For each monothetic subgroup $M=\overline{\<g\>}$
of $G$ there is an open normal subgroup $\Omega$ of $\Gamma$ which fixes
the elements of $M$ elementwise and the finite group
$\Gamma/\Omega$ acts on $M$ with the same orbits as $\Gamma$.

\item{\rm(ii)} The orbits  of $\Gamma$ on $\g=\L(G)$ for the
induced automorphic action of $\Gamma$ on $\L(G)$ are finite.

}
\msk
\bf
Proof.\quad\rm 
(i) Assume that $A=\overline{\<g\>}$ is a monothetic subgroup. If 
$\alpha\in \Gamma$, then $\alpha\in \Gamma_g$ means $\alpha\.g=g$,
that is, $g$ belongs to the fixed point subgroup Fix$(\alpha)$ of $\alpha$.
This is equivalent to $A\subseteq {\rm Fix}(\alpha)$, i.e.,
$\alpha\in \Gamma_a$ for all $a\in A$. So we have 
$\Gamma_g\subseteq \Gamma_a$ for all $a\in A$. Then the normal finite
index subgroup $\Omega_g=\bigcap_{\gamma\in\Gamma}\gamma\Gamma_g\gamma^{-1}$
is contained in all $\Gamma_a$, $a\in A$. This establishes (i) as
the remainder is clear.
\ssk
(ii) A compact abelian group $A$ is monothetic if and only if 
there is a morphism
$f\colon \Z\to A$ with dense image exactly when (in view of Pontryagin duality)
there is an injective morphism $\hat f\colon\hat A\to\hat\Z\cong \T$.
In the same spirit a compact abelian group is solenoidal
if and only if  there is a morphism
$f\colon \R\to A$ with dense image exactly when (in view of Pontryagin duality)
there is an injective morphism $\hat f\colon\hat A\to\hat\R\cong \R$. 
As $\hat A$ is discrete, this happens if and only if 
 $\hat A$ is algebraically a subgroup 
of $\R$. Now $\T$ has a subgroup algebraically isomorphic to $\R$
(see [\HM], Corollary A1.43). Thus if $\hat A$ can be homomorphically
injected into $\R$ it can be homomorphically injected into $\T$. 
Therefore
\ssk
\cen{\it every solenoidal compact group is monothetic.}
\ssk
Thus let $X\in\g$ be a one-parameter subgroup of $G$. Then by (i),
the image $X(\R)$ is contains in a monothetic subgroup. Hence by (i) above
there is an open normal subgroup $\Omega$ of $\Gamma$ such that each 
$\alpha\in\Omega$ satisfies $\alpha\.X(t)=X(t)$ for all $t\in\R$.
Hence $\alpha\.X=X$ in $\g$ with respect to the action of $\Gamma$ induced on
$\g=\L(G)$. Hence $\Gamma_X\supseteq\Omega$ 
has finite index for the action of $\Gamma$ on $\g$. \qed

Lemma 3.2(ii) completes the argument that the $G$-module $\g$ in the
case of a compact {\rm FC}-group $G$ satisfies $\fin \g=\g$.
\msk

In the following we refer to to the 
representation theory of compact groups as presented in [\HM], Chapters 3 
and 4.

\msk
\bf
Lemma 3.3.\quad \it Let $V$ be a locally convex weakly complete vector
space isomorphic to $\R^I$ and an effective
 $\Gamma$-module for a profinite group
$\Gamma$. In particular, the associated representation 
$\pi\colon \Gamma\to \Aut V$ is injective.
 Assume that $\fin V=V$. Then 

{\parindent2em

\item{\rm(i)} $V$ is a finite direct sum
(and product) $V_1\oplus\cdots \oplus V_k$ of isotypic components.

\item{\rm(ii)} $\Gamma$ is finite.

}
\msk
\bf
Proof. \quad \rm
(i) By  [\HM], Theorem 4.22,  $\fin V$ is a direct sum of its isotypic
components $V_\epsilon$, where $\epsilon\in \hat\Gamma$ is an equivalence 
class of irreducible representations of $\Gamma$. Each $V_\epsilon$ is a
module retract of $V$ under a canonical projection $P_\epsilon$
 and is therefore complete, and thus as a closed
vector subspace of a weakly complete vector space is weakly complete.
The universal property of the product 
$W\defi\prod_{\epsilon\in\hat\Gamma}V_\epsilon$ gives us an equivariant 
$\Gamma$-module morphism $\phi\colon V\to W$ of weakly complete
$\Gamma$-modules such that $\pr_\epsilon \ \circ \ \phi =P_\epsilon$. Since
morphisms of weakly complete vector spaces have a closed image
(see [\HMpro], Theorem A2.12(b)) and $$\sum_{\epsilon\in\hat\Gamma}V_\epsilon
=\{(v_\epsilon)_{\epsilon\in\hat\Gamma}\in W: v_\epsilon =0
\hbox{ for all but finitely many }\epsilon\in\hat\Gamma\}\subseteq W$$ 
is in the image of $\phi$ and is dense in
$W$ we know that $\phi$ is surjective. 
Now $V=\fin V$ and $\phi(\fin V)\subseteq \fin W$ whence $\fin W=W$.
It is readily verified that 

\cen{$W_\epsilon=\{(v_\eta)_{\eta\in\hat\Gamma}\in W:
v_\eta=0\hbox{ for }\eta\ne\epsilon\}$.}

\nin Therefore 
$\fin W=\sum_{\epsilon\in\hat\Gamma}V_\epsilon$, and thus
$$\sum_{\epsilon\in\hat\Gamma}V_\epsilon=W
=\prod_{\epsilon\in\hat\Gamma}V_\epsilon.$$
This equation, however, implies that the set
$S\defi\{\epsilon\in\hat\Gamma: V_\epsilon\ne\{0\}\}$ is finite. 
List the members of $\{V_\epsilon :\epsilon\in S\}$ 
as $V_1,\dots, V_k$. Then assertion (i) follows.
\ssk 
(ii) For a simple $\Gamma$-module $F$ denote by 
$[F]\in \hat\Gamma$
its equivalence class.  
Let $F_j$, $j=1,\dots,k$  be  simple $\Gamma$-modules 
such that $S=\{[F_j]: j=1,\dots,k\}$, $V_j=V_{[F_j]}$. 
 Then  let 
$\pi_j\colon \Gamma\to {\rm GL}(F_j)$, $j=1,\dots,k$
denote the associated representations
defined  by $\pi_j(g)(v)=g\.v$. 
By [\HM], Theorem 4.22,  $V_j\cong\Hom_\Gamma(F_j,V)\otimes F_j$, where
$g\in\Gamma$, $f\in \Hom(F_j,V)$, and $v\in V$ imply  $g\.(f\otimes v)
=f\otimes g\.v$. Since $V$ is an effective (or faithful) $\Gamma$-module, 
(i) implies 
$\bigcap_{j=1}^k\ker \pi_j=\{1\}$. 
However, a simple $\Gamma$-module is finite-dimensional (see
[\HM], Theorem 3.51), and so GL$(F_j)$ is a Lie group. 
Hence $\Gamma/\ker\pi_j\cong\im\pi_j$ is a compact 
profinite Lie group and is therefore finite.
Hence $\ker\pi_j$ is open in $\Gamma$. Thus $\{1\}$, being
a finite intersection of open subgroups, is open in $\Gamma$. 
Therefore, $\Gamma$ is discrete
and compact, hence finite.\qed

Thus we get the following result: 
\msk
\bf Lemma 3.4.\quad\it Let $G$ be a compact {\rm FC-}group. Then
the centralizer $Z(G_0,G)$ of the identity component is open.
\msk
\bf
Proof.\quad\rm  By Lemma 3.2 we can
apply Lemma 3.3 with $\Gamma=G/Z(G_0,G)$ and the
$\Gamma$-module $\g=\L(G)=\L(G_0)$. We conclude that $\Gamma$ is
finite. This proves the assertion.\qed
\msk

Let us explain the structure of a central extension of a compact
abelian group by a profinite group:

\msk
\bf Proposition 3.5. \quad\it 
Let $G$ be a compact group such that $G_0$ is
central. Then there is a profinite normal subgroup $\Delta$ such that
$G=G_0\Delta$. In particular, 
$$G\cong{G_0\times\Delta\over D},\quad D=\{(g,g^{-1}): g\in G_0\cap\Delta\}.$$
\msk\bf
Proof.\quad\rm See [\HM], Theorem 9.41 and Corollary 9.42.\qed

\msk
In order to pursue the structure of compact {\rm FC}-groups further
we cite Lemma 2.6, p. 1281 
from the paper by {\smc Shalev} [\shalev], proved with the aid of
the Baire Category Theorem:
\msk
\bf
Lemma 3.6.\quad \it If $G$ is a profinite {\rm FC}-group then
its commutator subgroup $G'$ is finite.\rm \qed
\msk 
In order to exploit this information, we need a further lemma:
\msk
\bf
Lemma 3.7.\quad \it If $N$ is a compact nilpotent group of class $\le2$ and
$N'$ is discrete, then the center $Z(N)$ has finite
index in $N$.
\ssk
\bf
Proof. \quad \rm Since  $N$ is nilpotent of class at most 2, we
have $N'\subseteq Z(N)$. Hence so  for each $y\in N$
the function $x\mapsto [x,y]:N\to N'$ is a morphism. Therefore,
if we set $A=N/N'$  we have a continuous $\Z$-bilinear map of
abelian groups $b\colon A\times A\to N'$, where
$b(xN',yN')= [x,y]$. Since $N'$ is discrete, 
$\{1\}$ is a neighborhood of the identity in $N'$.
On the other hand, $b(\{1_A\}\times A)=\{1\}$.
So for each $a\in A$ we have open neighborhoods $U_a$ of $1_A$
and $V_a$ of $a$, respectively, such that $b(U_a\times V_a)=\{1\}$
As $A$ is compact, we find a finite set $F\subseteq A$ of elements
such that $A=\bigcup_{a\in F}V_a$. Let $U=\bigcap_{a\in F} U_a$; then
$U$ is an identity neighborhood of $A$ and $b(U\times A)=\{1\}$.
Since $b$ is bilinear, this implies $b(\<U\>\times A)=\{1\}$. Then
the full inverse image $M$ of $\<U\>$ in $N$ under the quotient
morphism $N\to A$ is an open subgroup of $N$ satisfying
$[M, N]=\{1\}$ and is, therefore, central. Thus the center
$Z(N)$ of $N$ is open  and so has finite index in $N$.  \qed
\ssk
We remark that the proof of this lemma resembles that of 
Proposition 13.11, p.~574 of [\HMpro].
From the preceding two pieces of information we derive

\msk
\bf
Proposition 3.8.\quad \it Let $G$ be a compact group whose commutator
subgroup $G'$ is finite, and let
$Z(G',G)$ be the centralizer of $G'$ in $G$. Then the
center  of $Z(G',G)$ is a characteristic abelian open subgroup 
$A_G\defi Z(Z(G',G))$ of $G$.

\msk
\bf
Proof.\quad \rm Let $f\colon G\to \Aut(G')$ denote the morphism
defined by $f(g)(x)=gxg^{-1}$ for $x\in G'$. Since $G'$ is finite by
Lemma 3.6, so is $\Aut \ G'$. Hence $f(G)$ is finite as well and
so  $G/\ker f\cong f(G)$ is finite. If follows that the centralizer
$C\defi Z(G',G)$ of $G'$ in $G$,
being equal to $\ker f$, has finite index in $G$ and thus
is open. The center $Z(G')=G'\cap C$ of $G'$ is finite abelian,
and so the isomorphism $C/Z(G')\to G'C/G' \subseteq G/G'$
shows that the commutator subgroup $C'\subseteq Z(G')$
is a  finite subgroup of $C$. Since $G'\subseteq Z(Z(G',G),G)=
Z(C,G)$, the subgroup $C'\subseteq G'\cap C$ is central in $C$. 
Hence $C$
is nilpotent of class at most two with a finite commutator subgroup. 
Now Lemma 3.7 shows that $Z(C)$ is open in $C$, and since $C$
is open in $G$ we know that $Z(C)$, a characteristic subgroup of $G$,
is open in $G$. \qed
\msk
By Lemma 3.6, this applies to all profinite FC-groups.
Collecting our information on compact {\rm FC}-groups, we can
establish the following main structure theorem on the structure of compact
FC-groups:
\msk
\bf
Theorem 3.9. \quad \rm (The Structure of Compact FC-groups)
\it Let $G$ be a compact group. Then the following
two statements are equivalent:

{\parindent2em

\item{$(1)$} $G$ is an {\rm FC}-group. 

\item{$(2)$} $G/Z(G)$ is finite, that is, $G$ is center by finite.

\item{$(3)$} The commutator subgroup $G'$ of $G$ is finite.

}
\msk
\bf 
Proof.\quad \rm  (1)\implies(2): 
By Lemma 3.4, $Z(G_0,G)$ is open and so by
 Proposition 3.5 there is a profinite {\rm FC}-group $\Delta$
commuting elementwise with $G_0$ such that \hfill\break
$Z(G_0,G)=G_0\Delta$. 
By Lemma 3.6 and 
Proposition 3.8, $A_\Delta$ is an open characteristic abelian subgroup
of $\Delta$ whence  $A\defi G_0A_\Delta$ is an open characteristic abelian 
subgroup of $G_0\Delta=Z(G_0,G)$. Since $Z(G_0,G)$ is normal in $G$,
the  subgroup $A$ is normal in $G$. Since $A$ is open in $Z(G_0,G)$
and this latter group is open in $G$ by Lemma 3.4, $A$ is an
open normal subgroup of $G$.  Thus $G$ is an abelian by finite FC-group.
We claim that $G$ is therefore a center by finite group: Indeed,
let $f\colon G/A\to\Hom(A,A)$ be defined by $f(gA)(a)=gag^{-1}a^{-1}$.
For each coset $\gamma=gA$ the image of $f(\gamma)$ is finite since
$G$ is an FC-group, and so $Z(g,G)=\ker f(\gamma)$ has finite index.
Therefore $Z(G)=\bigcap_{\gamma\in G/A}\ker f(\gamma)$ has finite index,
and this proves the claim. 

\ssk

(2)\implies(1): Since $Z(G)\subseteq Z(g,G)$ for each $g$, by (2),  the
quotient $G/Z(g,G)$ is finite for all $g$. Hence $G$ is an FC-group. 

(1)\implies(3):
Let $G$ be an FC-group. By the equivalence of (1) and (2) we know
that $G_0$, being contained in the
open subgroup $Z(G)$, is central. By Proposition 3.5, there is a
profinite subgroup $\Delta$ such that $G=G_0\Delta$. Then
$[G,G]=[\Delta,\Delta]$. Now by Lemma 3.6, $[\Delta,\Delta]$ is
finite. 

(3)\implies(1): Set $C(g)=\{xgx^{-1}:x\in G\}$;
then $C(g)g^{-1}\subseteq G'$, whence $C(g)\in G'g$. Thus the finiteness of
$G'$ implies that of $C(g)$ for all $g\in G$. \qed
\msk  

In fact, a group $G$ satisfying the equivalent conditions of
 Theorem 3.9 above is what has been called
a BFC-{\it group}, that is, an FC-group with all conjugacy
classes of elements having bounded length.
\msk
Let us remark that center by 
finite groups are subject to classical central extension theory.
For instance,
if $Z(G)$ happens to be divisible by $q=|G/Z(G)|$, then $G=Z(G)E$ for
a finite subgroup $E$ such that $Z(G)\cap E$ is contained in 
$\{z\in Z(G): z^q=1\}$. 
(Cf.~method of proof of Theorem 6.10(i) of [\HM].)
\msk

At this point we return to the issue of 
probability of commuting elements
and thus to Borel probability measures on $G$.
Recalling Definition 2.1 we now  summarize our discussion as follows:

\msk
\bf 
 Theorem 3.10.\quad\it
Let $G$ be a compact group and $F$ its {\rm FC}-center.
Further
let $\mu_1$ and $\mu_2$ two  Borel probability measures on $G$
 and set $P=\mu_1\times\mu_2$
and $D\defi\{(g,h)\in G\times G:[g,h]=1\}$.
Assume that $\mu_j$ respects closed subgroups for $j=1,2$ and
that, if $G$ is not a Lie group,  $\mu_2$ respects even Borel subgroups.
Then the following 
conditions are equivalent:

{\parindent2em

\item{$(1)$} $P(D)>0$.

\item{$(2)$} $F$ is open in $G$. 

\item{$(3)$} The characteristic abelian subgroup $Z(F)$ is open in $G$.

}

Under these conditions, the centralizer $Z(F,G)$ of $F$ in $G$ is open,
and the finite group $\Gamma\defi G/Z(F,G)$ is finite and
acts effectively on $F$ with the same orbits as $G$ under the well defined
action  $\gamma\.x= gxg^{-1}$ for $(\gamma,x)\in \Gamma\times F$,
$g\in\gamma$.  The isotropy group
$\Gamma_x$ at $x\in  F$ is $Z(x,G)/Z(F,G)$, and the  set
$F_\gamma$ of fixed points under the action of $\gamma$ is $Z(g,F)$
for any $g\in \gamma$.
\msk
\bf
Proof.\quad\rm (1)$\Leftrightarrow$(2): 
This is a part of Corollary 3.1.

(2)$\Rightarrow$(3): 
The FC-center $F$ of $G$ is an  FC-group in
its own right. Then Theorem 3.9 shows that $Z(F)$ is open in $F$.
By (2), $F$ is open in $G$. Then
$Z(F)$ is open, and this establishes (3).

The implication (3)$\Rightarrow$(2) is trivial. 

Now assume that these conditions are satisfied. Then the open
subgroup $Z(F)$ is contained in the the centralizer $Z(F,G)=
\bigcap_{x\in F}Z(x,G)$, which is the kernel of the morphism $G\to \Aut F$
sending  $g\in G$ to $x\mapsto gxg^{-1}$. The homomorphism
$\pi\colon \Gamma=G/Z(F,G)\to \Aut F$ is therefore well-defined
by $\pi(\gamma)(x)=gxg^{-1}$, $\gamma=gZ(F,G)$, 
independently of the choice of
the representative $g\in \gamma$. If $x\in F$ and $\gamma\in \Gamma$,
say $\gamma=gZ(F,G)$, then $\gamma\in \Gamma_x$ iff $\gamma\.x=x$
iff $gxg^{-1}=x$ iff $g\in Z(x,G)$ regardless of the choice of 
$g\in\gamma$. Thus $\Gamma_\gamma= Z(x,G)/Z(F,G)$.
Similarly, we have $x\in  F_\gamma$ iff $\gamma\.x=x$ iff
$gxg^{-1}=x$ for $g\in\gamma$ iff $x\in Z(g,F)$. \qed
\msk

Later we shall  specialize the measures $\mu_1$ and $\mu_2$ in a reasonable way
and then search for useful measures $\mu_1$ and $\mu_2$ respecting
appropriate classes of subgroups.

\bsk 
Meanwhile, we shall specialize $\mu_1$ and $\mu_2$ both to Haar measure on $G$
and derive the main theorem on commuting elements in a compact group.

This case represents the special
class of examples in which $P(E)$ is
the probability that two randomly
picked elements commute. In this situation  one has called
$P(\{(x,y)\in G\times G: [x,y]=1\})$ the 
{\it commutativity degree} $d(G)$ of $G$.

Here, by the Main Theorem 3.10, the finite group
$\Gamma=G/Z(F,G)$ acts effectively on $F$ so that its orbits are
the $G$-conjugacy clases of elements of $F$ and that for the isotropy
groups and fixed point sets we have
$$\Gamma_x =Z(x,G)/Z(F,G)\quad\hbox{and}\quad (\forall g\in\gamma)\,
F_\gamma=Z(g,F).$$
In particular, if $\nu_F$ is Haar measure of $F$, for the closed 
subgroup $Z(g,F)$ of $F$ and $g\in\gamma\in\Gamma$ we conclude
$$\nu_F(F_\gamma)=\cases{0&if $F_\gamma$ has no inner points in $F$,\cr
       |F/Z(g,F)|^{-1}&if $F_\gamma$ is open in  $F$.\cr}$$

We now formulate and prove the following 
fundamental result:
\msk\bf
Main Theorem 3.11. \quad\it
Let  $G$ be any
 compact group and denote by  $d(G)$ its commutativity degree. 
Then we have the following conclusions:
\msk
{\smc Part}{\rm (i)} The following conditions are equivalent:
{\parindent2em

\item{$(1)$} $d(G)>0$.

\item{$(2)$} The center $Z(F)$ of the {\rm FC}-center $F$ of $G$
is open in $G$.} 

\msk
{\smc Part}{\rm(ii)} Assume that these conditions  of {\rm(i)} are satisfied. 
 Then there is a finite set of elements 
$g_1,\dots,g_n\in G$, $n\le|G/Z(F,G)|$,  such that
$$d(G)={1\over |G/F|\.|G/Z(F,G)|}\.\sum_{j=1}^n|F/Z(g_j,F)|^{-1}.$$
\msk
{\smc  Part}{\rm(iii)} $d(G)$ is always a rational number.

\msk
\bf
Proof.\quad\rm (i) follows directly from  3.10. 
For a proof of (ii) we let $\nu_G$ and $\nu_F$ be the Haar measures 
of $G$ and $F$, 
respectively,  and recall from 3.1 that 
$$ d(G)=P(D)=\int_{x\in F}\nu_G(Z(x,G))d\nu_G(x)    
            =\int_{x\in F}|G/Z(x,G)|^{-1}d\nu_G(x).$$
We note that the $\nu_G$-measure
of $F$ is $|G/F|^{-1}$; thus $\nu_F=|G/F|\.(\nu_G|F)$.
Hence
$$d(G)={1\over|G/F|}\.\int_F |G/Z(x,G)|^{-1}d\nu_F(x).$$
Now 
$$G/Z(x,G)\cong (G/Z(F,G))/(Z(x,G)/Z(F,G))=\Gamma/\Gamma_x,$$
and so, letting
$$E_F=\{(\gamma,x)\in\Gamma\times F: \gamma\.x=x\}$$
and $P=\nu_\Gamma\times\nu_F$, by $(*)$ above, we get 
$$d(G)={1\over|G/F|}\.\int_F|\Gamma/\Gamma_x|^{-1}d\nu_F
={1\over|G/F|}\. \int_F\nu_\Gamma(\Gamma_x) d\nu_F
={P(E_F)\over|G/F|}.$$
Next we select a set $\{g_1,\dots,g_n\}$ of elements of $G$ 
such that the cosets $g_jZ(F,G)$ are exactly those elements 
$\gamma\in\Gamma$  whose fixed point set 
$F_\gamma$ is open in $F$. Then $\Gamma_{g_jZ(F,G)}=Z(g_j,F)$ for all
$j=1,\dots,n$.
We also recall $\Gamma=G/Z(F,G)$ and  apply Lemma 2.5 to see that
$$P(E_F)={1\over|G/Z(F,G)|}\sum_{j=1}^n|F/Z(g_j,F)|^{-1}.$$
Therefore
$$d(G)={1\over |G/F|\.|G/Z(F,G)|}\.\sum_{j=1}^n|F/Z(g_j,F)|^{-1}.$$
This is what we claimed.
\msk
Part (iii) is now an immediate consequence of (ii) as $d(G)$ is trivially
rational if it is $0$.\qed\rm
 \msk
The main result says: {\it If the probablity that two randomly picked
elements commute in a compact group is positive, then, no matter how
small it is,  the group is almost abelian.}

We draw attention to the fact that in 3.11(ii) the commutativity
degree of $G$ is expressed in purely arithmetic terms via the 
group theoretical data  $F$, $Z(F,G)$, and $Z(g_j,F)$.
In the much simpler case of finite groups  this aspect was discussed in
[\lesc]. 

\bsk

\cen{\smc 4. Examples}
\msk
Finally, 
we record some examples. Recall 
that 
$D(G)=\{(x,y)\in G\times G: [x,y]=1\}$,
that  $P=\nu\times\nu$ with normalized Haar measure $\nu$  on $G$,
and that 
$d(G)\defi P(D(G))$ is the  commutativity degree.
All of the examples in our list are metabelian, and the first two
are in fact nilpotent. The Main Theorem 3.11 explains why this is not too far 
from the most general situation.

The first example is trivial and arises from finite groups:

\msk
\bf Example 4.1.\rm  \quad Let $H$ be the
8-element quaternion group and $P=\nu\times\nu$. Then $d(H)=5/8$.
Let $A$ be any compact connected group, e.g., the circle group $\T$.
Let $G=A\times H$. The  random selection of pairs of commuting elements in
$A$ and the random selection of pairs of commuting elements in $H$ 
are independent events, and thus 
$$d(G)= P(D(A))\times P(D(H))= d(A)\.d(H)= 5/8.$$
In view of Theorem 3.9, this example is a compact FC-group.
\msk

\bf
Example 4.2.\rm\quad
 Let $H$ be the class 2 nilpotent compact
group of all $3\times3$-matrices
$$M(a,b;z)\defi\pmatrix{1 & a & z \cr
           0 & 1 & b \cr
           0 & 0 & 1 \cr},$$
where $a, b, z$ range through the ring $\Z_p$ of $p$-adic integers.
Let $Z$ be the closed central subgroup of $H$ of all $M(0,0;pz)$ ,
$z\in\Z_p$. Then $G=H/Z$ is a compact nilpotent $p$-group of class 2
whose commutator group $[H,H]/Z$ contains all elements $M(0,0;z)Z$
and thus is isomorphic to $\Z(p)=\Z/p\Z$. Its center $Z(G)$ consists of
all elements $M(a,b;z)Z$ with $a, b\in p\Z_p$, $z\in\Z_p$, 
whence $G/Z(G)\cong\Z(p)^2$. 
The factor group $G/[G,G]$ is isomorphic to $\Z_p^2$. The
subgroup of all $M(a,0;0)$, $a\in \Z_p$, is isomorphic to $\Z_p$; it
is topologically generated by $M(1,0;0)$. By Theorem 3.9 again,
this example yields a compact FC-group which is not the product
of a central group with a finite group. \qed

{\smc Erfanian} and the second author have shown the following
proposition, whose proof draws on 
ideas of {\smc Lescot} [\lesc] on finite groups.

\msk
\bf
Proposition 4.3. \rm ([\ER], Theorem A)\quad {\it Assume $G$ is an
nonabelian compact {\rm FC}-group  such that  $G/Z(G)\cong (\Z/p\Z)^2$ 
for a prime $p$. Then $d(G)=(p^2+p-1)/p^3.$}\ssk

 The group $G$ in Example 4.2 therefore
has the commutativity degree $(p^2+p-1)/p^3$. 
All  examples of compact groups $G$
in the existing literature having commutativity degree
 $d(G)=(p^2+p-1)/p^3$ arise
from direct products of finite groups such as e.g.\  
Example 4.1 above.
In this sense  $G=M(a,b;z)$ exhibits a first ``nontrival'' explicit
compact example illustrating Proposition 4.3.

Similar formulae as in 4.3 can be found  in [\EK, \ER, \gust].

\msk
The following example is a  nonnilpotent metabelian group.
\msk
\bf
Example 4.4.\rm\quad 
Let $A$ be a compact abelian group and set $A_a=\{a\in A:2\.a=0\}$.
 Define $G=A\sdir{}\{1,-1\}$ with multiplication
$$(a,\epsilon)(a',\epsilon')=(a+\epsilon\. a',\epsilon\epsilon').$$
Then $Z((a,1),G)=A\times\{1\}$ for $2a\ne0$, and
$Z((0,-1),G)=A_2\times\{-1,1\}$. The orbits are
$G\.(t,1)=\{(\pm a,1)\}$,  $G\.(0,-1)=2\.A\times\{-1\}$.
Then $$\eqalign{%
D(G)&= (A\times\{1\})^2\cup (A\times\{-1\})\times(A_2\times\{1\})\cr
&\cup (A_2\times\{1\})\times(A\times\{-1\})\cup (A_2\times\{-1\})^2,\cr}$$ 
and so, setting $t=\nu_A(A_2)$, we get 
$$d(G)={1\over4}+{1\over2}t+{1\over4}t^2=\left({1+t\over2}\right)^2.$$
\msk
If $A_2$ is without inner points in $A$ we have $t=0$ and so $d(G)={1\over4}$.
For $A=\T$ we get the ``continuous dihedral group'' which shows, 
among other things,
that under the circumstances of Theorem 3.10, 
we may have $Z(F)=F\ne G$ and that $G$ may not be
center by finite. 

Note that the power function $g\mapsto g^2$ 
is constant equal to 1.

\msk

\bsk \cen{\smc 5. On the history and background of the problem} 
\msk
A study of the probability that two randomly picked elements $x$ and
$y$ of a compact group $G$ commute was initiated by W.~H.~Gustafson
in [\gust]. Thereby he extended to the infinite case an idea which
was put forward a few years earlier by P.~Erd\H os and P.~Tur\'an in
the context of finite groups (s.\ [\gust]  again and [\ert]).
Perhaps Erd\H os and Tur\'an attempted to reformulate in a
statistical way a famous problem, posed by Paul~Erd\H os himself
and solved  partially by Bernhard Neumann in [\neum]:

{\it For any class of groups $\cal{X}$, let $\cal{X}^*$ denote the class
of all groups $G$ such that every infinite subset of $G$ contains a
pair of distinct elements which generate an $\cal{X}$-subgroup of
$G$. Is it true  $\cal{X}^* \subseteq \cal{X}$?}

In [\neum], Bernhard Neumann answered this problem in the negative
by showing that for the class $\cal{A}$ of all abelian groups, the
class $\cal{A}^*$ turns out to be  the class of finite central
extensions of abelian groups. His methods dealt with combinatorial
techniques, coverings of suitable sets and remarks on the size of
 centralizers. Exactly the same methods can be found in [\ert,
\gust].
 The approach of P.~X.~Gallagher for
computing  the
probability that two randomly picked elements of a finite group 
commute is completely different in so far as he 
used character theory (s. [\lesc]) to express
this probability  in terms of the number of
conjugates of a given element. More details can be found in [\diac]
and [\lesc]. In particular, Diaconis' survey [\diac] 
illustrates the fact that considerable
information on the structure of groups  may be obtained along this route,
and it helped to
 motivate the growing interest in these issues that has emerged in the
literature in recent years. Regarding infinite groups,
compact topological groups were the natural  ones to which the results on
 finite groups could be extended 
 as soon as appropriate methods suitable for the topological
context were found. We have  seen that the concepts of measure theory serve
as a substitute for  an approach through  character theory. However
no visible contributions  so far deal with this point, even though 
the  literature on the representation theory of compact
groups is considerable.

In [\HM], one approaches the subject of this paper from the other end,
following up on a classical theme. In [\HM], Proposition 6.86
and Exercise E6,18 is is shown that in a nonabelian connected compact
Lie group $G$ the set of all pairs $(x,y)\in G\times G$ such that
$\<x,y\>$ is free nonabelian and dense in $G$ is (i) dense in $G\times G$,
(ii) has a meager complement, and (iii) has Haar measure 1 in $G\times G$.
Thus one has no chance that the commutativity degree 
of such groups is positive,
as is amply confirmed through the principal results of this article.

Peter M. Neumann [\neumtwo] studied the set 
$D(G)=\{(x,y)\in G\times G: [x,y]=1\}$
in a finite group and consequences arising from the assumption that
the commutatiity degree $d(G)=P(D(G))$ is positive in this case.
In [\levai], Levai and Pyber extended his results to profinite, that is,
 totally disconnected compact groups 
and  proved results similar to the ones in the present paper.  

The papers [\EK, \ERE, \ER] illustrate that it is possible to extend
W.~H.~Gustaf\-son's initial idea  to more sophisticated notions of
probability on compact groups. In principle, for any finite sequence
of words
$w_\alpha(x_1,x_2,\ldots,x_n)$, $\alpha=1,\dots,N$ in $n$ 
variables and for any compact group $G$
the probability of the set
$$\{(x_1, x_2,\ldots, x_n)\in G^n:  w_1(x_1,x_2,\ldots,x_n)=
\cdots=w_N(x_1,\dots,x_n)=1\}$$
in the group $G^n$ becomes a reasonable object of study. 
For instance, in [\EK] and
[\ER] the authors consider the words 
$w_{ij}(x_1,\ldots,x_n)=[x_i,x_j]$ for
$1\leq i<j\le n$ and the measure of the set
$$\{(x_1,\ldots,x_n)\in G^n:\ [x_i,x_j]=1, 1\leq i< j\le n\}.$$
Some simple observations concerning the relation $w=1$ on compact Lie groups
are made in [\HM], Lemma 6.83.
\msk

{\bf Acknowledgments.} We are grateful  to Karl-Hermann
Neeb and Georg W.~Hofmann  for posing inspiring  questions, to
Ben Green for pointing out and clarifying reference [\neumtwo], 
leading to the paper [\levai] by L. Levai and L. Pyber
which had been brought
 to our attention  through a lecture by Aleksander Ivanov;
while they treat profinite groups only, in the area 
of positive probability for commuting pairs there is 
overlap between their results and ours. We thank
 an anonymous referee of an earlier version
of this paper for detecting and communicating
flaws which we eliminated in the present one.
\bsk 

{\bf Literature}

\ssk
[\bour] Bourbaki, N., {\it Int\'egration}, Chap.~7 et 8,
Hermann, Paris,  1963.

\ssk
[\bred] Bredon, G., {\it Introduction to Compact Transformation
Groups}, Academic Press, New York, 1972.

\ssk
[\diac] Diaconis, P., Random walks on groups: characters and
geometry, in: {\it Groups St. Andrews 2001 in Oxford}, Vol. I,
London Math. Soc. Lecture Note Ser. {\bf304}, Cambridge Univ. Press,
Cambridge, 2003, pp.~120--142.

\ssk [\ert] Erd\H os, P., and P. T\'uran, On some problems of
statistical group theory, {\it Acta Math. Acad. Sci. Hung.} {\bf 19}
(1968), 413--435.

\ssk
[\EK] Erfanian, A.,  and R.~Kamyabi--Gol,   On the mutually
commuting $n$-tuples in compact groups, {\it Int.\ J.\ Algebra} 
{\bf1} (2007), 251--262.

\ssk
[\ERE] Erfanian, A.,  and  R.~Rezaei, On the commutativity
degree of compact groups, {\it Arch.\ Math.\ (Basel)} {\bf 93} (2009),
201--212.

\ssk
[\ER] Erfanian, A.,  and  F.~Russo, Probability of mutually
commuting $n$-tuples in some classes of compact groups, {\it Bull.\ 
Iran.\ Math.\ Soc.} {\bf 34} (2008), 27--37.

\ssk
[\gust]  Gustafson, W.\ H., What is the probability that two
group elements commute? {\it Amer.\ Math.\ Monthly} {\bf 80} (1973),
1031--1304.

\ssk
[\hew] Hewitt, E., and K. Ross, 
Abstract Harmonic Analysis, Vol.I, Springer, Berlin, 1963.

\ssk
[\HM] Hofmann, K. H., and S.~A.~Morris, {\it The Structure of
Compact Groups}, de Gruyter, Berlin, Second Edition 2006.

\ssk
[\HMpro]  ---, The Lie Theory of connected Pro-Lie Groups,
                Eur.\ Math.\ Soc.\ Publ.\ House, 2007.

\ssk
[\hofmosii] ---,
Elements of Compact Semigroups, Charles E. Merill, Columbus,
Ohio, 1966, xii+384 pp.

\ssk
[\lesc] Lescot, P., Isoclinism classes and commutativity
degrees of finite groups, {\it J.\ Algebra} {\bf 177} (1985),
847--869.

\ssk
[\levai] L\'evai, L., and L. Pyber, Profinite groups with many
commuting pairs or involutions, {\it Arch. Math.} (Basel) 
{\bf75} (2000), 1--7.

\ssk
[\neum] Neumann, B. H., On a problem of P.~Erd\H os, {\it J.\
Aust.\ Math.\ Soc.\ Ser.\ A} {\bf 21} (1976), 467--472.

\ssk
[\neumtwo] Neumann, P. M., Two combinatorial problems in group theory.
{\it Bull. London Math. Soc.} {\bf21} (1989), 
456-458.

\ssk
[\shalev]
Shalev, A., {\it Profinite groups with restricted centralizers},
Proc.\ Amer.\ Math.\ Soc.\ {\bf122} (1994), 1279--1284.

\ssk
[\tD]   tom Dieck, T., {\it Transformation Groups}, 
de Gruyter, Berlin, 1987.

\vskip1cm

Authors' addresses:
\bsk 
\obeylines
\line{%
\vbox{\hsize.52\hsize
Karl H. Hofmann, 
Fachbereich Mathematik
Technische Universit\"at
Schlossgartenstr. 7
64289 Darmstadt, Germany
hofmann@mathematik.tu-darmstadt.de}
\hfill
\vbox{\hsize.45\hsize
Francesco G. Russo
Department of Mathematics
University of Palermo
via Archirafi 34
90123 Palermo, Italy
francescog.russo@yahoo.com}}

\end